\documentstyle [12pt]{article}
\input{epsf.sty}
\begin{document}

\title{First passage percolation and a model for competing spatial growth}
\author{Olle H\"{a}ggstr\"{o}m\thanks{Research supported by grants from the
Swedish Natural Science Council and the Royal Swedish Academy of Sciences.}
\and Robin Pemantle\thanks{Research supported in part by a grant from the
Alfred P. Sloan Foundation, by a Presidential Faculty Fellowship and by NSF
grant \# DMS9300191} \and \\
Chalmers University of Technology \\
University of Wisconsin}
\date{}
\maketitle

\begin{abstract}
An interacting particle system modelling competing growth on the ${\bf Z}^2$ 
lattice is defined
as follows. Each $x\in{\bf Z}^2$ is in one of the states $\{0,1,2\}$. 
$1$'s and $2$'s remain in their states forever, while a $0$ flips to a
$1$ (resp.\ a $2$) at a rate equal to the number of its neighbours which
are in state $1$ (resp.\ $2$). This is a generalization of the
well known Richardson model. $1$'s and $2$'s may be thought of as two
types of infection, and $0$'s as uninfected sites. We prove that if we 
start with a single site in state $1$ and a single site in state $2$,
then there is positive probability for the event both types of infection
reach infinitely many sites. This result implies that the spanning tree
of time-minimizing paths from the origin in first passage percolation
with exponential passage times has at least two topological ends with
positive probability. 


\end{abstract}

\newpage

\section{Introduction and statements of results}

We consider two-dimensional first passage percolation, 
where each pair of nearest neighbours of ${\bf Z}^2$ have an edge
connecting them and each edge is equipped with a nonnegative random variable
which is interpreted as the time it takes to traverse the edge (see \cite{D}
and \cite{K} for reviews). We specialize to the case where these passage
times are i.i.d.\
exponentials; this is the so called Richardson model \cite{R}. By time
scaling, we may assume without loss of genertality that the exponential
distribution has mean one. 
For $x,y\in{\bf Z}^2$, let $T(x,y)$ be the time taken to reach $y$ from $x$,
i.e.\ $T(x,y)$ is the infimum over all paths from $x$ to $y$ of
the sum of the passage times along the path. Write ${\bf 0}$ for the origin
and let, for $t\geq 0$,
\[
B(t)=\{x\in{\bf Z}^2: T({\bf 0},x)\leq t\}
\]
so that $B(t)$ is the set of sites reached from the origin by time $t$. 
It follows from the
memoryless property of the exponential distribution that $B(t)$ is a Markov
process, and the following interacting particle system formulation is
natural. Each site $x$ takes the value $\eta_x(t)\in\{0,1\}$ at time $t$. Let 
\[
\eta_x(t)=\left\{
\begin{array}{ll}
1 & \mbox{if } x\in B(t) \\
0 & \mbox{otherwise.}
\end{array} \right.
\]
and write $\eta(t)=\{\eta_x(t)\}_{x\in{\bf Z}^2}$ for the whole configuration 
of 0's and 1's at time $t$. 
We may think of sites in state 0 as healthy and those in state 1 as infected. 
Each site in state 1 remains in this state forever and tries to infect each of 
its four nearest neighbours at rate one, so that a site in state 0 flips
to state 1 at rate equal to the number of nearest neighbours with state 1. 
At time 0, only the origin is infected. 

One of the most important results in first passage percolation is the 
existence of the so called time constant $\mu$ such that
\begin{equation} \label{Kingman}
\lim_{n\rightarrow\infty}\frac{T({\bf 0},(n,0))}{n}=\mu 
\mbox{\hspace{4 mm} a.s.}
\end{equation}
which is a consequence of Kingman's subadditive ergodic theorem. This result
in fact holds under much more general conditions than those considered in 
this paper. The precise value of $\mu$ is not known, although some bounds are
available. We will need a lower bound on $\mu$, and the best such bound we
are aware of is
\begin{equation} \label{Janson}
\mu>0.298
\end{equation}
as shown in \cite{J}. 
Similar time constants exist for all directions. Moreover, there 
is the following 
asymptotic shape result which, somewhat loosely speaking, states that
analogues of (\ref{Kingman}) hold for all directions simultaneously: Let 
\[ 
\bar{B}(t)=\{x\in {\bf R}^2: \exists y \in[\textstyle{-\frac{1}{2},
\frac{1}{2}}]^2, z\in B(t)
\mbox{ such that } y+z=x\},
\]
i.e.\ $\bar{B}(t)$ is a ``fattened''
version of $B(t)$. Then there exists a nonrandom compact convex set $B_0$
such that for all $\varepsilon>0$
\[
(1-\varepsilon)B_0\subseteq \frac{\bar{B}(t)}{t} \subseteq (1+\varepsilon)B_0
\mbox{\hspace{2 mm} eventually a.s.,}
\]
see e.g.\ \cite{D}.

We will be interested in a different aspect of the evolution of infection, 
namely the {\bf tree of infection}, to be denoted by $\Gamma$. Let $\Gamma(t)$ 
be the graph with vertex set $B(t)$ and edge set obtained as follows. For
each $x\in B(t)\setminus \{{\bf 0}\}$, let $e_x$ be the edge connecting $x$
to the vertex $y$ from which $x$ got infected, and let 
$\{e_x: x \in B(t)\setminus \{{\bf 0}\}\}$ be 
the edge set of $\Gamma(t)$. Each vertex of $\Gamma(t)$ gets a unique
path to ${\bf 0}$, so $\Gamma(t)$ is 
indeed a tree. Moreover, both the vertex set and
the edge set of $\Gamma(t)$ are increasing in $t$, so the limiting object
$\Gamma=\lim_{n\rightarrow\infty}\Gamma(t)$ also exists. The tree structure
of $\Gamma(t)$ is inherited by $\Gamma$. Equivalently, we may define $\Gamma$ 
as the graph with vertex set ${\bf Z}^2$ and edge set given by
\[
\bigcup_{x\in {\bf Z}^2\setminus \{{\bf 0}\}}\{e: e \mbox{ is an edge of the 
fastest path from ${\bf 0}$ to }x\}.
\]
Our main interest is in the number of topological ends of $\Gamma$, 
i.e.\ how many infinite self-avoiding paths starting at ${\bf 0}$
does $\Gamma$ contain? Let $K(\Gamma)$ denote the number of such paths.
By a standard compactness argument, $K(\Gamma)\geq 1$. Newman \cite{N} 
has shown that $K(\Gamma)=\infty$ a.s.\ provided a certain
hypothesis concerning uniformly bounded curvature of $B_0$. (See also 
\cite{LN} for related results.)
The uniformly bounded curvature hypothesis
is highly plausible, but has so far not been not proved. We shall prove
the following result which is a small step towards the conjecture that
$K(\Gamma)=\infty$ a.s.

\medskip\noindent
{\bf Theorem 1.1:} {\sl The number of topological ends $K(\Gamma)$ of the tree
of infection satisfies
\[
P(K(\Gamma)\geq 2) \geq {\textstyle \frac{4\mu-1}{3}}>0.
\]}

In order to prove this result, we shall first study a simple and natural
model for competing spatial growth. We will now go on to describe this
model, which is a variant of the Richardson model and which we propose
to call the {\bf two-type Richardson model}. 

Consider an interacting particle system on ${\bf Z}^2$ with state space 
$\{0,1,2\}$, where 0's may be thought of as healthy sites, and 1's and 2's
as two different types of infection. The evolution is as follows. A site in
state 1 stays in this state forever, and the same thing holds for a site
in state 2. Both types of infected sites try to infect each of their nearest
neighbours at rate one, so that a site in state 0 flips to state 1 at rate
equal to the number of nearest neighbours in state 1, and to state 2
at rate equal to the number of nearest neighbours in state 2. 
We will write $\xi(t)$ for the configuration at time $t$; $\xi(t)$ will be
an element of $\{0,1,2\}^{{\bf Z}^2}$. Note that
if we disregard the type of infection (i.e.\ if we watch this system
evolve wearing a pair of glasses which prevents us from distinguishing 
between 1's and 2's), then the system behaves exactly as the ordinary
(one-type) Richardson model. This is an immediate consequence of the fact 
that the total flip rate of a site in state 0 equals its total number of 
infected nearest neighbours.

We start at time 0 with an infection of type 1 at {\bf 0} and an 
infection of type 2 at another site $x$, all other sites being healthy. 
The model can then be described in first passage percolation terms:
A site $y$ is infected at time $T(\{{\bf 0},x\},y)$, which we define as the
infimum, over all paths starting at ${\bf 0}$ or $x$ and ending at $y$, of
the sum of the passage times along the path. Since the distribution of
the passage times of the edges is continuous, it is not hard to see that
the infimum is in fact a.s.\ a minimum which is attained
for a unique path. If this fastest path starts at ${\bf 0}$, then $y$
gets infection of type 1, otherwise it gets type 2. 

We may think of the two-type Richardson model as a crude model for
two growing bacteria colonies (or two political empires) competing for space. 
It may happen that at some early stage, one of the types of
infection completely surrounds the other type which then is prevented from
growing indefinitely. Write $A$ for the event that this does not happen,
in which case both types of infection will grow indefinitely. The first 
question one would like to answer 
about the two-type Richardson model is whether or not $P(A)>0$
(it is obvious that $P(A)<1$).
The answer to this question is in fact independent of $x$, as stated in
the following proposition, 
where we write $P_{{\bf 0},x}(A)$ for the probability
of $A$ with the described starting configuration. The proof will be given in
Section 2. 

\medskip\noindent
{\bf Proposition 1.2:}
{\sl For any $x_1,x_2\in{\bf Z}^2\setminus {\bf \{0\}}$, 
\[
P_{{\bf 0},x_1}(A)>0 \mbox{\hspace{4 mm} implies \hspace{4 mm}}
P_{{\bf 0},x_2}(A)>0.
\] }

\medskip\noindent
Hence we may restrict attention to $x={\bf 1}:=(1,0)$, i.e.\
to the two-type Richardson model starting with a 1 at the origin, a 2 right
next to the origin, 
and all other sites healthy. Our main result on the two-type
Richardson model now says

\medskip\noindent
{\bf Theorem 1.3:}
\[
P_{{\bf 0},{\bf 1}}(A)\geq {\textstyle\frac{4\mu-1}{3}} >0,
\]

\medskip\noindent
so that with positive probability both bacteria colonies
keep growing forever. Inserting (\ref{Janson}) 
yields $P_{{\bf 0},{\bf 1}}(A)>0.064$. 

Let us now go back to the one-type Richardson model and the question of the
number $K(\Gamma)$ of topological ends of the tree of infection. Start with
a single infected site at ${\bf 0}$ and suppose that the first site that
${\bf 0}$ infects is ${\bf 1}$. It is then an immediate consequence of
Theorem 1.3 and the identification 
between the one-type Richardson model and the
two-type Richardson model with types disregarded, that there is positive
probability for the event that $\Gamma$ has two different self-avoiding
paths to infinity: one which goes through ${\bf 1}$ and one which does not.
Hence {\bf Theorem 1.3 implies Theorem 1.1}. 

Similarly, if it had turned out that infinite growth of both bacteria 
colonies had probability zero, then we would have been able to conclude that 
$K(\Gamma)=1$ a.s. This remark is of course an empty statement, but we point
it out anyway because it is feasible that a similar implication might
be useful in some of the possible extensions of the model considered here;
see below. 

In order to prove Theorem 1.3, we will need the following proposition, which
will be proved in Section 2. We think it is of some interest in
its own, and it also seems related to questions concerning the roughness of
the boundary of $B(t)$, studied e.g.\ in \cite{K2} and \cite{NP}. 

\medskip\noindent
{\bf Proposition 1.4:}
{\sl For any $\varepsilon>0$, there exist infinitely many 
$x=(x_1,x_2)\in{\bf Z}^2$ in the right half-plane, such that
\begin{equation}\label{PreProp}
P[T({\bf 0},(x_1,x_2))>T({\bf 0},(x_1-1,x_2))] > 
{\textstyle \frac{2\mu+1}{3}} - \varepsilon,
\end{equation}
which is greater than $0.5$ when $\varepsilon$ is small. 
In fact, the stronger result holds that for some $l\in\{0,1,\ldots\}$
\begin{equation} \label{Prop}
\limsup_{n\rightarrow\infty}P[T({\bf 0},(n,l))>T({\bf 0},(n-1,l))] >
{\textstyle \frac{2\mu+1}{3}}-\varepsilon.
\end{equation} 
}

\medskip\noindent
Inserting (\ref{Janson}) gives that the right hand sides in
(\ref{PreProp}) and (\ref{Prop})
can be made greater that $0.532$. Intuitively, Proposition 1.4 says that
there are sites in ${\bf Z}^2$ arbitrarily far away from the origin which 
``strongly feel'' from which
direction the infection is coming. It seems reasonable to believe that 
the $\limsup$ in (\ref{Prop}) in fact is a limit, and independent of $l$.   

Apart from proving the conjecture about $\Gamma$ a.s.\ having infinitely
many topological ends, 
there are various other ways in which one might want to extend and improve the 
results of this paper. An obvious question is what happens when ${\bf Z}^2$
is replaced by ${\bf Z}^d$ for $d\geq 3$. Another direction is to allow
passage time distributions other than the exponential. The Markovian
behaviour of $B(t)$ is then lost whence the interacting particle system
formulation becomes less natural, but the questions about $\Gamma$ remain
just as natural as for the exponential case. An inspection of the proof of
Proposition 1.4 shows that the result $P(K(\Gamma)\geq 2)>0$ extends to the
case where the passage times of the edges have a continuous distribution
with a hazard rate which is bounded between two constants $\alpha<\beta$
satisfying $\frac{\alpha}{\beta}>\frac{1-3\mu}{\mu}$, (here $\mu$ still
denotes the time constant when the passage time distribution is
a mean one exponential). One can also ask what happens if we extend
the two-type Richardson model to having three or more different types,
in the obvious way. For this extension, showing that $k$
different types simultaneously can grow indefinitely from a finite
starting configuration is equivalent to showing that
$P(K(\Gamma)\geq k)>0$. Yet another direction of generalization 
would be to allow the two types to have different infection rates
$\lambda_1$ and $\lambda_2$, and for this extension we conjecture that
$P(A)=0$ whenever $\lambda_1\neq \lambda_2$. 

Since Theorem 1.3 implies Theorem 1.1, 
it only remains to prove Propositions 1.2 and 1.4, and Theorem 1.3.
The remainder of the paper is devoted to this task. The hard work is in
the proof of Proposition 1.4.  

\section{Proofs}

{\bf Proof of Proposition 1.2:}
Suppose $P_{{\bf 0},x_1}(A)>0$, and let $S_r$ be the circle of radius
\[
r=\max(|x_1|,|x_2|)+2
\]
(here $|\cdot |$ is the Euclidean norm) 
centered at ${\bf 0}$.
By conditioning on the first $n$ infections (for some $n$ which we need not
define explicitly here)
we can find $t>0$ and a configuration 
$\xi\in\{0,1,2\}^{{\bf Z}^2}$ such that $\xi$ 
restricted to ${\bf Z}^d\cap S_r$ contains only 1's and 2's, and such that
$P_{{\bf 0},x_1}(\xi(t)=\xi,A)>0$. 
Note that $\{x\in{\bf Z}^2: \xi_x=1\}$ and $\{x\in{\bf Z}^2: \xi_x=2\}$
both must be connected sets (in the graph-theoretic sense with edges between
all nearest neighbours).
Let 
\[
B_\xi=\{x\in{\bf Z}^2: \xi_x\in\{1,2\}\}
\]
be the set the set of infected sites in $\xi$; by choice of $r$ also $B_\xi$
is connected. Define 
\[
\partial B_\xi^\ast=\{x\in{\bf Z}^2: \xi_x\in\{1,2\} \mbox{ and }
\exists y\in{\bf Z}^2 \mbox{ such that } \xi_y=0\mbox{ and }|x-y|=1 \},
\] 
i.e.\ $\partial B_\xi$ is 
the set of infected sites with at least one healthy neighbour. It is easy to
check, using the strong Markov property, that $P_{{\bf 0},x_1}(A|\xi(t)=\xi)$
does not depend on the values of $\xi$ on $B_\xi\setminus \partial B_\xi$. 
By choice of $r$, we may now construct a
configuration $\xi^\prime\in\{0,1,2\}^{{\bf Z}^2}$ such that
\begin{description}
\item{(i)}
$B_{\xi^\prime}=B_\xi$ (whence in particular $\partial B_{\xi^\prime}=
\partial B_\xi$),
\item{(ii)}
$\xi^\prime=\xi$ on $\partial B_\xi$, 
\item{(iii)}
$\xi^\prime_{\bf 0}=1$ and $\xi^\prime_{x_2}=2$,
\item{(iv)}
$\{x\in{\bf Z}^2: \xi^\prime_x=1\}$ and $\{x\in{\bf Z}^2: \xi^\prime_x=2\}$
are both connected.
\end{description}
$B_\xi^\prime$ is finite, whence by (iii) and (iv) we have $P_{{\bf 0},x_2}
(\xi(t)=\xi^\prime)>0$ because only finitely many infections are involved. 
Hence, 
\begin{eqnarray*} 
P_{{\bf 0},x_2}(A) & \geq & P_{{\bf 0},x_2}(\xi(t)=\xi^\prime, A) \\
& = & P_{{\bf 0},x_2}(A|\xi(t)=\xi^\prime)P_{{\bf 0},x_2}(\xi(t)=\xi^\prime) \\
& = & P_{{\bf 0},x_1}(A|\xi(t)=\xi)P_{{\bf 0},x_2}(\xi(t)=\xi^\prime) \\
& > & 0
\end{eqnarray*}
as desired. $\Box$

\medskip\noindent
{\bf Proof of Proposition 1.4:}
Our aim is to show that (\ref{Prop}) holds for fixed $\varepsilon>0$ 
and some $l$.  
We will consider the evolution of $B(t)$ inside the half-strip
$S_k=\{0,1,\ldots\}\times\{0,1,\ldots,k\}$ for some large $k$ 
(taking $k>\frac{2}{3\varepsilon}$ will suffice). Writing $N_k(t)$
for the number of infected sites in $S_k$ at time $t$, we have, as an easy
consequence of (\ref{Kingman}), that
\begin{equation} \label{lim4}
\lim_{t\rightarrow\infty}\frac{N_k(t)}{t}=
{\textstyle \frac{k+1}{\mu}} \mbox{\hspace{4 mm} a.s.}
\end{equation}
For $x\in{\bf Z}^2$, write $\lambda_x(t)$ for the flip rate of $x$ at time $t$,
and write $\lambda_{S_k}(t)$ for $\sum_{x\in S_k}\lambda_x(t)$. It follows from
(\ref{lim4}) and an application of the strong law of large numbers that 
\begin{equation} \label{lim5}
\lim_{t\rightarrow\infty}{\textstyle\frac{1}{t}}\int_0^t
\lambda_{S_k}(u)du={\textstyle \frac{k+1}{\mu}} \mbox{\hspace{4 mm} a.s.}
\end{equation}
Now write $N^G_k(t)$ ($G$ as in good) for the number of infected sites in
$S_k$ at time $t$ which became infected after their nearest neighbour to the
left, and write $N_k^B(t)$ ($B$ for bad) for the corresponding number of
sites which became infected after their nearest neighbour to the right. 
By keeping track of the left-right nearest neighbour pairs, we have, as for
(\ref{lim4}), 
\[
\lim_{t\rightarrow\infty}\frac{N^G_k(t)+N^B_k(t)}{t}
={\textstyle \frac{k+1}{\mu}} \mbox{\hspace{4 mm} a.s.}
\]
although note that $N^G_k(t)+N^B_k(t)=N_k(t)$ 
holds only approximately because an 
individual site may sometimes be infected after neither or both of its
left-right nearest neighbours. Write $\lambda^G_{S_k}(t)$ for 
$\sum\lambda_x(t)$ where this time the sum is taken over 
uninfected sites in $S_k$ whose
nearest neighbour to the left is already infected at time $t$, and define 
$\lambda^B_{S_k}(t)$ analogously. We have as in (\ref{lim5}) that
\[
\lim_{t\rightarrow\infty}\frac{\int_0^t
\lambda^G_{S_k}(u)du}{N_k^G(t)}=1 \mbox{\hspace{4 mm} a.s.}
\]
and
\begin{equation} \label{lim6}
\lim_{t\rightarrow\infty}\frac{\int_0^t
\lambda^G_{S_k}(u)du+\int_0^t
\lambda^B_{S_k}(u)du}{\int_0^t
\lambda_{S_k}(u)du}=1 \mbox{\hspace{4 mm} a.s.}
\end{equation}
In order to show that (\ref{Prop}) holds for some $l\in\{0,\ldots,k\}$, 
it suffices to prove that
\[
\liminf_{t\rightarrow\infty}
\frac{N^G_k(t)}{N_k(t)}>
{\textstyle\frac{2\mu+1}{3}}-\varepsilon \mbox{\hspace{4 mm} a.s.}
\]
To this end, it is sufficient to show that
\begin{equation} \label{lim7}
\liminf_{t\rightarrow\infty}
\frac{\int_0^t\lambda^G_{S_k}(u)du}{\int_0^t\lambda_{S_k}(u)du}>
{\textstyle\frac{2\mu+1}{3}}-\varepsilon \mbox{\hspace{4 mm} a.s.} 
\end{equation}
so this is what we will proceed to prove.

We now look at a configuration $\tilde{\eta}_{S_k}\in\{0,1\}^{S_k}$, which will
serve as an example of what $\eta_{S_k}(t)$ may look like at some reasonably
large time point $t$. Since
$\eta_{(0,i)}(t)=1$ for $i=0,\ldots,k$ and all sufficiently large $t$ a.s.,
we assume that $\tilde{\eta}_{(0,i)}=1$ for $i=0,\ldots,k$. Also 
$\eta_x(t)=0$ for fixed $t$ and all but finitely many $x\in{\bf Z}^2$, so we
furthermore assume that $\tilde{\eta}_x=0$ for all but finitely many 
$x\in S_k$. Now pick $i\in\{0,1,\ldots,k-1\}$ and write $\tilde{\eta}_{S_k,i}$
for $\tilde{\eta}_{S_k}$ restricted to $\{0,1,\ldots\}\times\{i,i+1\}$,
i.e.\ $\tilde{\eta}_{S_k,i}$ is the configuration on two infinite
adjacent horizontal rows. Given $\tilde{\eta}_{S_k,i}$, we may partition
$\{0,1,\ldots\}$ into 1-blocks, 0-blocks, down-blocks, up-blocks and
excursion-blocks as in Figure 1. 

\begin{figure}[hptp]
\epsfxsize=7cm
\hfill{\epsfbox{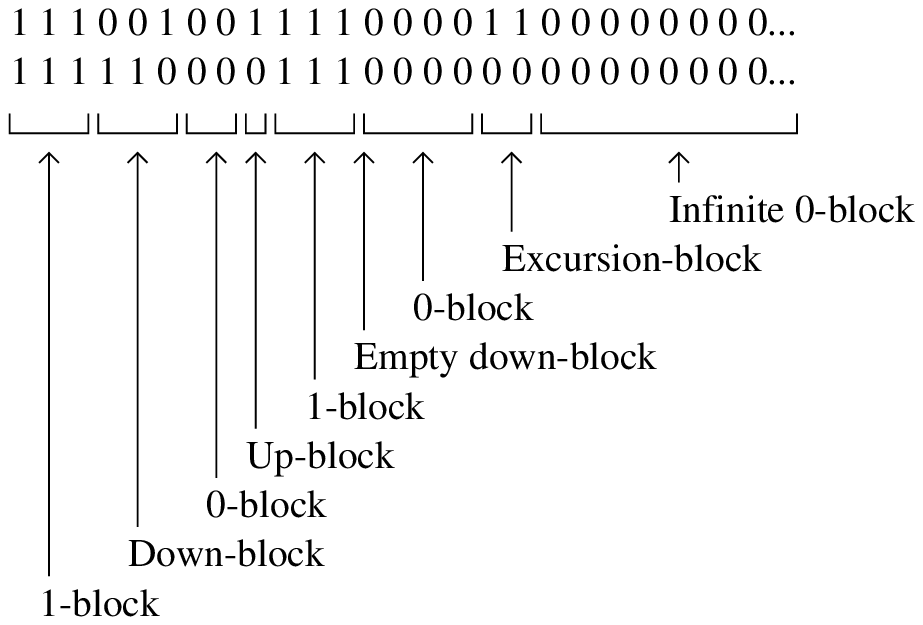}\hfill}
\caption{A typical blockpartition.}
\end{figure}

The blocks are determined (uniquely, given $\tilde{\eta}_{S_k,i}$) as follows.
A 1-block is a maximal sequence of $[\mbox{\small$1\atop 1$}]$'s and a
0-block is a maximal sequence of $[\mbox{\small$0\atop 0$}]$'s, while
down-blocks, up-blocks and excursion-blocks are maximal sequences consisting of
$[\mbox{\small$0\atop 1$}]$'s and $[\mbox{\small$1\atop 0$}]$'s. All blocks
are finite except for the final 0-block. A 1-block is always followed either by
an up-block or by an excursion-block. A finite 0-block is always followed by
an up-block or an excursion-block. A down-block is always followed by a 0-block
while an up-block always is followed by a 1-block. An excursion-block, finally,
is always surrounded either by two 1-blocks or by two 0-blocks. Down-blocks
and up-blocks (but no other blocks) may be empty; an empty down-block is 
inserted whenever a $[\mbox{\small$1\atop 1$}]$ is followed by a
$[\mbox{\small$0\atop 0$}]$, and, similarly, an empty up-block is inserted
whenever a $[\mbox{\small$0\atop 0$}]$ is followed by a
$[\mbox{\small$1\atop 1$}]$. 

\noindent
Call 1-blocks and 0-blocks {\bf pure}, and call other blocks {\bf mixed}. 
Call a pair of nearest neighbours {\bf hot} if one of the sites has value 1
and the other has value 0. Each hot pair will be associated with a mixed
block as follows. If a $[\mbox{\small$1\atop 1$}]$ is followed by a
$[\mbox{\small$0\atop 0$}]$, then the two hot pairs are associated with the 
corresponding empty down-block, and similarly when a 
$[\mbox{\small$0\atop 0$}]$ is followed by a $[\mbox{\small$1\atop 1$}]$.
Any other hot pair intersects exactly one mixed block and is then associated
with this block. 

We now consider rates of infection between sites in $\{0,1,\ldots\}\times
\{i,i+1\}$ when $\eta_{S_k}(t)=\tilde{\eta}_{S_k}$. Write $\tilde{\lambda}(M)$
for the weighted total infection rate for the hot pairs associated with the
mixed block $M$, with weight $\frac{1}{2}$ for horizontal infections and
weight 1 for vertical infections. The reason for these weights is that 
each horizontal infection is accounted for in $\tilde{\eta}_{S_k,i}$ for
two different values of $i$ (this is not true for rows $0$ and $k$, but that
will be corrected for later). Write $\tilde{\lambda}^G(M)$ (resp.\
$\tilde{\lambda}^B(M)$) for the corresponding sum where only good 
(resp.\ bad) infections are counted. The following is easily 
checked to hold for $\tilde{\lambda}(M)$, $\tilde{\lambda}^G(M)$ and
$\tilde{\lambda}^B(M)$. 
\[
\begin{array}{lrr}
\mbox{Type of block} & \tilde{\lambda}(M) & 
\hspace{4 mm}\tilde{\lambda}^G(M)-\tilde{\lambda}^B(M) \\ \\
\mbox{empty down-block} & 1 & 1\hspace{14 mm} \\
\mbox{empty up-block} & 1 & -1\hspace{14 mm} \\
\mbox{nonempty down-block} & \geq 2 & 2\hspace{14 mm} \\
\mbox{nonempty up-block} & \geq 2 & -2\hspace{14 mm} \\
\mbox{excursion-block} & \geq 2 & 0\hspace{14 mm}
\end{array}
\]
Write $\tilde{\lambda}_{S_k,i}$ (resp.\ $\tilde{\lambda}^G_{S_k,i}$
and $\tilde{\lambda}^B_{S_k,i}$) for $\tilde{\lambda(M)}$ (resp.\ 
$\tilde{\lambda}^G(M)$ and $\tilde{\lambda}^B(M)$) summed over all mixed
blocks in $\tilde{\eta}_{S_k,i}$. Note that since $\tilde{\eta}_{S_k,i}(t)$
starts with a 1-block and ends with a 0-block, we have that the number of
down-blocks exceeds the number of up-blocks by exactly 1. Using this
observation, and
the above table, we may check (e.g.\ via induction over the
number of down-blocks) that
\begin{equation} \label{lim8}
\tilde{\lambda}_{S_k,i}\geq 4-3(\tilde{\lambda}^G_{S_k,i}-
\tilde{\lambda}^B_{S_k,i})
\end{equation}
with equality if and only if the only mixed blocks in $\tilde{\lambda}_{S_k,i}$
are empty down-blocks and nonempty up-blocks with $\tilde{\lambda}(M)=2$ for
each up-block $M$. 
Applying (\ref{lim8}) to $\eta_{S_k,i}(t)$ and summing over $i$, we get
\[
\sum_{i=0}^{k-1}\lambda_{S_k,i}(t)\geq 4k-3\sum_{i=0}^{k-1}
(\lambda^G_{S_k,i}(t)-\lambda^G_{S_k,i}(t))
\]
for all sufficiently large $t$ a.s., whence
\begin{equation} \label{lim9}
\liminf_{t\rightarrow\infty}\sum_{i=0}^{k-1}{\textstyle\frac{1}{t}}\int_0^t
(\lambda_{S_k,i}(u)+3\lambda^G_{S_k,i}(u)-3\lambda^B_{S_k,i}(u))du\geq 4k.
\end{equation}
Since the number of infections up to time $t$ in rows $0$ and $k$ is
asymptotically $2\mu^{-1}t$, we have that
\[
\liminf_{t\rightarrow\infty}{\textstyle\frac{1}{t}}\int_0^t
\left( \sum_{i=0}^{k-1}\lambda_{S_k,i}(u)-\lambda_{S_k}(u) \right)du \leq
2\mu^{-1}
\]
and for the same reason
\[
\liminf_{t\rightarrow\infty}{\textstyle\frac{1}{t}}\int_0^t
\left( \sum_{i=0}^{k-1}\lambda^B_{S_k,i}(u)-\lambda^B_{S_k}(u) \right)du \leq
2\mu^{-1}.
\]
Hence, we may modify (\ref{lim9}) to get
\[
\liminf_{t\rightarrow\infty}{\textstyle\frac{1}{t}}\int_0^t
(\lambda_{S_k}(u)+3\lambda^G_{S_k}(u)-3\lambda^B_{S_k}(u))du
\geq 4k-6\mu^{-1}.
\]
Combining this with (\ref{lim6}) yields 
\[
\liminf_{t\rightarrow\infty}{\textstyle\frac{1}{t}}\int_0^t
(6\lambda^G_{S_k}(u)-2\lambda_{S_k}(u))du
\geq 4k-6\mu^{-1}
\]
i.e.
\begin{eqnarray*}
\liminf_{t\rightarrow\infty}{\textstyle\frac{1}{t}}\int_0^t
6\lambda^G_{S_k}(u)du & \geq &
4k-6\mu^{-1}+2\lim_{t\rightarrow\infty}{\textstyle\frac{1}{t}}\int_0^t
\lambda_{S_k}(u)du \\
& = & 4k-6\mu^{-1}+2(k+1)\mu^{-1} \\
& = & k(4+2\mu^{-1})-4\mu^{-1}
\end{eqnarray*}
which in conjunction with (\ref{lim5}) implies
\[
\liminf_{t\rightarrow\infty}\frac{\int_0^t\lambda^G_{S_k}(u)du}
{\int_0^t\lambda_{S_k}(u)du}\geq {\textstyle 
\frac{k(4+2\mu^{-1})-4\mu^{-1}}{6k\mu^{-1}}
=\frac{2\mu+1}{3}-\frac{2}{3k}}.
\]
Taking $k>\frac{2}{3\varepsilon}$,
this implies (\ref{lim7}), so the proof is complete. $\Box$

\medskip\noindent
{\bf Proof of Theorem 1.3:}  
Suppose for contradiction that 
$P_{{\bf 0},{\bf 1}}(A)< \frac{4\mu-1}{3}-2\varepsilon$ for some 
$\varepsilon>0$. 
Then, by symmetry,
the probability that infection of type 2 eventually stops growing is
greater than 
$\frac{1}{2}(1-\frac{4\mu-1}{3}+2\varepsilon)=\frac{2-2\mu}{3}+\varepsilon$. 
This implies
\[
\limsup_{n\rightarrow\infty}P[T({\bf 0},(n,l))>T({\bf 1},(n,l))]< 
{\textstyle 1-\frac{2-2\mu}{3}}-\varepsilon
={\textstyle\frac{2\mu+1}{3}}-\varepsilon,
\]
where $l$ can be chosen as in Proposition 1.4. 
By reflecting the realization of passage times 
in the line $x_1=\frac{n}{2}$, we see
that the pair
$(T({\bf 0},(n,l)),T({\bf 1},(n,l)))$ has the
same joint distribution as
$(T({\bf 0},(n,l)),T({\bf 0},(n-1,l)))$, so that
\[
\limsup_{n\rightarrow\infty}P[T({\bf 0},(n,l))>T({\bf 0},(n-1,l))] <
{\textstyle \frac{2\mu+1}{3}} - \varepsilon
\]
contradicting (\ref{Prop}). $\Box$


\begin{thebibliography}{9}

\bibitem{D} Durrett, R. (1988) {\sl Lecture Notes on Particle Systems and
Percolation}, Wadsworth \& Brooks/Cole, Pacific Grove. 

\bibitem{J} Janson, S. (1981) An upper bound for the velocity of
first-passage percolation, {\sl J. Appl. Probab.} {\bf 18}, 256--262.

\bibitem{K} Kesten, H. (1987) Percolation theory and first-passage percolation,
{\sl Ann. Probab.} {\bf 15}, 1231--1271.

\bibitem{K2} Kesten, H. (1993) On the speed of convergence in
first-passage percolation, {\sl Ann. Appl. Probab.} {\bf 3}, 296--338. 

\bibitem{LN} Licea, C. and Newman, C.M. (1996) Geodesics in two-dimensional
first-passage percolation, {\sl Ann. Probab.} {\bf 24}, 399--410. 

\bibitem{N} Newman, C.M. (1995) A surface view of first-passage percolation,
{\sl Proceedings of the 1994 International Congress of Mathematicians} (ed.\
S.D. Chatterij), 1017--1023,
Birkh\"auser, Boston.

\bibitem{NP} Newman, C.M. and Piza, M. (1995) Divergence of shape fluctuations
in two dimensions, {\sl Ann. Probab.} {\bf 23}, 97--1005. 

\bibitem{R} Richardson, D. (1973) Random growth in a tesselation, {\sl Proc.
Cambridge Phil. Soc.} {\bf 74}, 515--528. 

\end{thebibliography}
\end{document}